\documentclass[11pt,a4paper,oneside,english]{article}

\usepackage[top=2cm, bottom=2cm, left=2cm, right=1.5cm]{geometry}

\usepackage[utf8]{inputenc}
\usepackage[T1]{fontenc}
\usepackage[english]{babel}

\usepackage{authblk}

\usepackage{amsmath}
\usepackage{amssymb}
\usepackage{amsthm}

\usepackage{mathtools}

\usepackage{algorithm}
\usepackage{algpseudocode}

\usepackage[capitalise,noabbrev]{cleveref}
\usepackage{url}

\usepackage{tikz}
\usepackage{tikz-qtree}
\usepackage{pgfplots}
\usepackage{pgfplotstable,booktabs}
\usetikzlibrary{calc,math}

\theoremstyle{plain}
\newtheorem{Theorem}{Theorem}[section]
\theoremstyle{definition}

\newtheorem{Remark}[Theorem]{Remark}
\newtheorem{Problem}[Theorem]{Problem}

\newcommand{\du}{w}
\newcommand{\dphi}{\psi}

\newcommand{\ddu}{\delta u}
\newcommand{\ddphi}{\delta \varphi}
\newcommand{\ddU}{\delta U}

\newcommand{\strain}{e}

\newcommand{\sprod}[2]{\left( #1 , #2 \right)}
\newcommand{\norm}[1]{\Vert #1 \Vert}

\usepackage{xparse}

\NewDocumentCommand{\summary}{m O{-1} O{-1}}{		
	\pgfplotstabletypeset[
	columns={lvl, dofs, cells, h, eps},
	columns/lvl/.style={column name=$\ell$},
	columns/dofs/.style={sci, zerofill, precision = 1},
	columns/h/.style={column name=$h$, sci, precision = 1},
	columns/eps/.style={column name=$\varepsilon$, sci, precision = 1},
	columns/cells/.style={column name={elements}, sci, zerofill, precision = 1},
	every head row/.style={before row=\toprule,after row=\midrule},
	every last row/.style={after row=\bottomrule},
	skip rows between index={#2}{#3},
	]{#1}
}

\newcommand{\MSC}[1]{\small \vspace{4mm} {\bf MSC:} #1}

\begin{document}
	
	\selectlanguage{english}
	\pagenumbering{arabic}
	
      \title{Parallel matrix-free higher-order finite element solvers\\ for phase-field fracture problems}
	
	\author[1]{D. Jodlbauer}
	\author[1]{U. Langer}
	\author[2,3]{T. Wick}
	
	\affil[1]{Johann Radon Institute for Computational and Applied Mathematics, Austrian Academy of Sciences, Altenbergerstr. 69, A-4040 Linz, Austria}
	\affil[2]{Institut f\"ur Angewandte Mathematik, Leibniz Universit\"at Hannover, Welfengarten 1, 30167 Hannover, Germany}
	\affil[3]{Cluster of Excellence PhoenixD (Photonics, Optics, and Engineering - Innovation Across Disciplines), Leibniz Universit\"at Hannover, Germany}
	
	\date{}
	
	\maketitle
	
	\begin{abstract}
          
        Phase-field fracture models lead to variational problems that
        can be written as
         a coupled variational equality and inequality system.
        Numerically, such problems can be treated with Galerkin finite elements and primal-dual active set methods.
        Specifically, low-order and high-order finite elements may be employed, where, for the latter, only
        few studies exist to date.
		The most time-consuming part in the discrete version of the primal-dual active set (semi-smooth Newton) algorithm 
		consists in the solutions of changing linear systems arising at each semi-smooth Newton step. 
		We propose a new parallel matrix-free monolithic multigrid preconditioner
		for these systems. We provide two numerical tests, and discuss the performance of 
		the parallel solver proposed in the paper.
		Furthermore, we compare our new preconditioner with a block-AMG preconditioner available in the literature.


		\MSC{2010 74R10; 65M60; 49M15; 35Q74}

	\end{abstract}

\section{Introduction}

	Many applications require the solution of partial differential equations (PDEs).
	Frequently, these are solved by Finite Element Methods (FEM), or related approaches like Isogeometric Analysis, 
	which discretize the continuous PDE.
	To obtain accurate solutions, 
	{we eventually need}
	to solve huge linear systems of equations.
	This 
	{may become}
	a challenging task
	{since}
	the computational effort increases rapidly.
	Hence, solvers that are able to handle large-scale linear systems are required.
	This gives rise to specialized solvers, that are adapted towards a specific PDE.
	In this work, we focus on 
	{efficient parallel solvers for problems in 
	{phase-field fracture (PFF) propagation.}

	The origin of crack modeling dates back to 
	{Griffith}
	\cite{Gr21} in 1921, who proposed the first model for fractures in brittle materials.
	The phase-field approach we are using was developed by Francfort and Marigo \cite{FrMa98}, who put the fracture model in an energy minimization context.
	Therein, the originally low-dimension crack is approximated by means of a continuous phase-field function.
	The 
	{PFF}
	model was extended in \cite{MiWeHo10,AmMaMa09,AmGeDe15} by decomposing the stress terms into tensile and compressive parts.
	{However, the}
	physically correct way to do this splitting
	is disputed within the fracture community.
	Examples exist where either one or the other method is superior.
	There has been {developed} a multitude of further extensions
	{of}
	the original model, e.g., {to} pressure-driven fractures \cite{MiWhWi15a,WiLa16,YoBo16} in porous-media \cite{MiMaTe15,MiWhWi19,HeMa17,ChBoYo19}, multi-physics \cite{MiMa16} and many more;
	see also the survey papers \cite{BoFr19,WuNgNg19,WhWiLe20} and 
        the monograph \cite{Wi20}.

	PFF problems have the advantage, that they reduce to a system of PDEs,  which can be solved by adapting well-known strategies.
	On the other hand, a smooth approximation of an originally sharp fracture comes with some drawbacks, e.g., crack-length computation or interface boundary conditions along the fracture.
	There exist different approaches than PFF that allow treatment of sharp cracks, e.g., \cite{BeBl99,MeDu07,AlJoVa11}.

	The numerical solution of PFF problems is particularly challenging, mainly due to the so-called no-heal condition.
	This ensures that a crack does not regenerate itself.
	Mathematically, this leads to an additional variational inequality.
	The solution gets further complicated by the non-convexity of the underlying energy functionals.
	Many different approaches have been proposed in the literature to handle PFF; see, e.g.,
	\cite{BoFrMa00,WhWiWo14,MiHoWe10,HeWhWi15,KoKr20}. 
	In this work, we use the primal-dual active-set method, a 
	{version}
	of the semi-smooth Newton method \cite{HiItKu03}, which was applied to PFF in \cite{HeWhWi15}. 

	The most time consuming part in the simulation of PDEs is the solution of the arising linear systems of equations.
	In the context of nonlinear PDEs, these appear after linearization within Newton's algorithm, or similar linearization approaches.
	Hence, multiple linear systems need to be solved per simulation step.
	For increasing problem sizes, this becomes a severe bottleneck, and well optimized solvers are required.
	In \cite{JoLaWi19}, we presented a solver based on the matrix-free framework of the C++ FEM library deal.II \cite{AlArBa18}.
	Further strategies treating the linear systems are presented in \cite{FaMa17}.
	Matrix-free methods avoid storing the huge linear system, which can become a real issue in terms of memory consumption.
	The linear systems are solved using the Generalized Minimum Residual (GMRES) iterative solver using a geometric multigrid preconditioner.
	This setup does not need explicit knowledge about the matrix entries of the linear system, but only need to perform matrix-vector multiplications.
	These can be carried out without actually storing the matrix.

	Matrix-free approaches are particularly favorable for high-polynomial shape-functions.
	In these cases, the number of non-zero entries within the matrix increases, rendering the storage of the matrix and sparse-matrix-vector multiplication (SpMV) a lot more expensive.
	On the other hand, matrix-free methods do not suffer from an increase polynomial degree.
	These approaches originate from the field of spectral methods \cite{Or80,DeFiMu02}, where high-order ansatz functions are used.
	
	One key ingredient to keep the computational time low is parallelization.
	In this work, we show the effects of different ways of parallelization.
	First, we are interested in the distributed solution of the whole problem.
	Therein, the computational work is split across multiple CPUs (ranks), each with its own independent memory.
	Necessary data exchange among CPUs is done using the Message Passing Interface (MPI).
	Using multiple CPUs gives us the possibility to fit larger problems into memory, as each core only stores small, almost independent parts of computational domain.
	All cores run in parallel, resulting in faster computations.
	Secondly, we explicitly utilize vector instruction sets provided by modern CPUs.
	These are special instructions for the CPU that are capable of performing multiple computations at once.
	This is more thoroughly described in \cref{sec:simd}.

	In this paper, we 
	{study}
	the performance of the AMG-based solver by \cite{HeWi18} and the matrix-free geometric multigrid from our previous work \cite{JoLaWi19}.
	Both approaches as well as the treatment of the nonlinearities are given in \cref{sec:solutionmethods}.
	In \cref{sec:compare}, we compare both approaches 
	{with respect to (wrt)}
	computational time, memory consumption and parallel efficiency.
	Details on the parallelization are presented in \cref{sec:parallel}.
	We focus on dependence on the polynomial degree of the finite elements.
	For studies regarding $h$-dependence, we refer to the respective papers.
	\Cref{sec:numericresults} describes the numerical examples and compares the results to those found in the literature to validate our implementations.
	Details on the phase-field model are briefly described in \cref{sec:problem}, for more information we refer to the corresponding literature.
	We also added a discussion on the derivative of the eigensystem in \cref{subsec:DerivativeEigensystem}.

\section{Preliminaries}
\label{sec:problem}

In the following, we will repeat the basic notation used in  \cite{JoLaWi19} as well as the governing equations. 
Furthermore, we will introduce the primal-dual active-set method that can be interpreted as 
semi-smooth Newton method.
Finally, we present the computation of the derivative of the eigensystem that plays a fundamental role in Miehe-type splittings.

\subsection{Phase-Field Fracture Model}
\label{subsec:PhaseFieldFractureModel}

	We will use $(a,b) \coloneqq (a,b)_{D} \coloneqq \int_{D} a \cdot b \ dx$ for the standard $L^2$ scalar product, and $(A,B) \coloneqq (A,B)_{D} \coloneqq \int_{D} A : B \ dx$ for tensor-valued functions.
	
	\Cref{fig:notation} shows the computational domain $D \subset \mathbb{R}^d$ with an inscribed fracture $\mathcal{C} \subset \mathbb{R}^{d-1}$.
	The remaining undamaged parts of $D$ are denoted by $\Omega$.
	We use the phase-field approach from \cite{BoFrMa00,AmTo90, AmTo92}, where the lower-dimensional fracture is approximated by using an additional 
	{scalar function}
	$\varphi : D \rightarrow [0,1]$.
	Values of $0$ denote a completely fractured domain, whereas $1$ represents undamaged parts.
	The size of the fracture approximation is controlled by the length-scale parameter $\varepsilon$.
	
	\begin{figure}[H]
		\centering
		\begin{tikzpicture}
	\draw (0,0) rectangle (4,4);
	
	\draw[
		rounded corners=20.0, dotted,
		fill={rgb:red,1;green,1;blue,1;white,100}
	] (0.4,1.2) rectangle (3.6,2.8);
	\draw[
		rounded corners=8.0, densely dotted, 
		fill={rgb:red,1;green,1;blue,1;white,15},
	] (1.0,1.7) rectangle (3.0,2.3);
	
	\draw[thick] (1.25,2) -- (2.75,2);
	
	\node at (2.9,2.0) {\tiny $\mathcal{C}$};
	\node at (2,2.15) {\tiny $\varphi \approx 0.0$};
	\node at (2,2.55) {\tiny $\varphi \approx 0.5$};
	\node at (2,3.3) {\tiny $\varphi = 1$};
\end{tikzpicture}
		\caption{
			Computational domain $D$ with inscribed fracture $\mathcal{C}$ and its phase-field approximation.			
		}		
		\label{fig:notation}
	\end{figure}
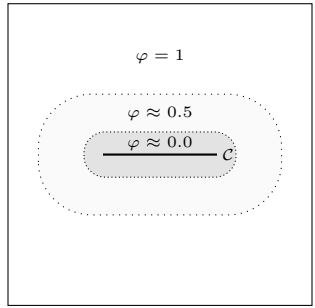
	
	The unknown displacement field $u: D \rightarrow \mathbb{R}^d$ and fracture $\varphi: D \rightarrow [0,1]$ are computed by means of an energy minimization problem
	\begin{equation*}
		\mathcal{E}(u, \varphi) = \int_{D} g(\varphi) E_s(\strain) dx + \frac{G_c}{2} (\frac{1}{\varepsilon} \norm{1 - \varphi}^2 + \varepsilon \norm{\nabla \varphi}^2) \rightarrow \min,
	\end{equation*}
	with additional displacement boundary conditions on $\partial D$.
	
	The decay of the material stiffness is controlled by the degradation function $g(\varphi)$.
	In this work, we use $g(\varphi) \coloneqq (1 - \kappa)\varphi^2 + \kappa$, with a small regularization parameter $\kappa \ll 1$ to avoid singularities in fractured domains.
	Different choices of $g$ can be found in the literature; see, e.g., \cite{KuScMu15,BoFrMa00,SaKeBe18,GeDe19}.
	
	The solid energy $E_s(\strain)$ is given by $\frac{1}{2} \lambda \mbox{tr}^2(\strain) + \mu (\strain,\strain)$ with the strain tensor $\strain \coloneqq \frac{1}{2}(\nabla u + \nabla u^T)$ and Lam\'e parameters $\mu$ and $\lambda$.
	
	We will further consider splitting the elastic energy into tensile and compressive parts $E_s = E_s^+ + E_s^-$ based on the spectral decomposition proposed by Miehe and coworkers in \cite{MiWeHo10}.
	The modified elastic energy functional then reads $E_s^+ \coloneqq \frac{1}{2} \lambda \langle \mbox{tr}(\strain) \rangle_+^2 + \mu \mbox{tr}(e_+^2)$, see \cite{MiWeHo10,MiHoWe10}.
	Here, $a^+ \coloneqq \max(0, a)$ for real numbers $a$, and $e^+ \coloneqq \sum_{i=1}^{d} \lambda_i^+ v_i \otimes v_i,$ for matrices $e$ with eigensystem $(\lambda_i, v_i)$.

	Consequently, the stress tensor splitting is defined as $\sigma^+(u) = \partial_e E_s^+$, i.e. $\sigma^+(e(u)) = \lambda \mbox{tr}(e)^+ I + 2 \mu e^+$.
	Various splitting schemes are discussed in \cite{BoVeSc12}[Section 2.2] and \cite{AmGeDe15}.
	
	This leads to the extended energy functional 
	\begin{equation*}
		\mathcal{E}(u, \varphi) = \int_{D} g(\varphi) E_s^+(\strain) + E_s^- dx + \frac{G_c}{2} (\frac{1}{\varepsilon} \norm{1 - \varphi}^2 + \varepsilon \norm{\nabla \varphi}^2) \rightarrow \min,
		\label{eq:energy_split}
	\end{equation*}
	
	Furthermore, fractures are only allowed to grow, i.e. $\partial_t \varphi < 0$.
	In our quasi-static formulation, this means that, for each loading step $n$, we need to enforce $\varphi^n \leq \varphi^{n-1} \quad \dot{\forall} x \in D$.
	The first-order necessary condition (Euler-Lagrange system) of \cref{eq:energy_split} is given by the variational inequality problem
	\begin{Problem}[Euler-Lagrange System]\label{eq:euler_lagrange}
		Find displacement $u \in V $ and fracture $\varphi \in W_{in}$ such that
		\begin{align*}
			& \left(  g (\varphi) \sigma^+(u), \strain(\du) \right) + \left(\sigma^-(u), \strain(\du) \right) = 0 \quad & \forall \du \in V, \\
			& \left(  \partial_\varphi g(\varphi) E_s^+(\strain(u)), \dphi - \varphi \right) & \\
			& + G_c \left( -\frac{1}{\varepsilon} ( 1-\varphi, \dphi - \varphi ) + \varepsilon \left( \nabla \varphi, \nabla (\dphi - \varphi ) \right) \right) \geq 0 \quad & \forall \dphi \in W \cap L^\infty(D),
		\end{align*}
		with spaces $V \coloneqq H^1_0(D)^d, W \coloneqq H^1(D), W_{in} \coloneqq \{ w \in W: w \leq \varphi^{old} \text{ a.e.} \}$; 
		see e.g. \cite{MiWhWi19}.
	\end{Problem}

	\begin{Remark} 
		In the literature, the phase-field approximation described here is also referred to as AT~2 model (named after Ambrosio/Tortorelli) introduced in \cite{AmTo90} for the Mumford-Shah problem and the original variational fracture formulation \cite{BoFrMa00}.
		Changing the last line in \cref{eq:euler_lagrange} leads to another well-known model referred to as AT~1; see, e.g., \cite{BoMaMa14}.
		However, in this work, we only consider the AT~2 model.
	\end{Remark}

\subsection{Treatment of the Variational Inequality Constraints}

	Solving such variational inequalities is a challenging task that admits many different solution techniques.
	Popular choices are penalization approaches, where one incorporates violations to the no-heal condition into the variational form, i.e. adding terms like $\gamma (\varphi - \varphi^{n-1})^+$ for large enough values of $\gamma$.
	Improvements as in \cite{WhWiWo14,Wi17} aim to determine $\gamma$ adaptively, or provide lower bounds as in \cite{GeDe19}.
	In \cite{MiHoWe10}, the inequality is replaced by a so-called history field $\mathcal{H}(x, t) \coloneqq \max \{ E_s^+(\strain(u)) \}$, which again leads to an equality system to be solved.

	In this work, we use the primal-dual active-set method for fractures as presented in \cite{HeWhWi15}.
	In the next section, we continue 
	{with the description of the} 
	solution approaches used in our simulations, including details 
	{to the finite element version of the
	active-set algorithm 
	mentioned above.}

\section{Discretization and Solution Algorithms}
\label{sec:solutionmethods}

	In this section we briefly discuss the discretization methods used.
	We start with a very short description of FEM.
	More details can be found in one of the standard textbooks.
	We continue with the key parts of our solution strategy: the primal-dual active-set method and the linear solvers used therein.
	In \cref{subsec:DerivativeEigensystem}, we illustrate how the derivative of the eigenvector and eigenvalues is computed in the case of Miehe-type splitting.

	\subsection{Finite Element Discretization}
	\label{sec:discretization}

	In particular, we use the Finite Element Methods on quadrilateral/hexahedral meshes to obtain discrete spaces for displacement and phase-field variables.
	The shape functions are chosen to be globally continuous piecewise polynomials of degree $p$ in every coordinate direction ($Q_p$).
	The code is implemented in C++ with the help of the FEM library deal.II \cite{AlArBa18}.
	We make particular use of the matrix-free geometric multigrid framework included, which is described in more detail in \cite{Ka04,JaKa11}.

	\subsection{Primal-Dual Active-Set Algorithm}
	\label{sec:activeset}

	The idea of the active-set algorithm is to enforce the constraints on a set $\mathcal{A}$, and perform a Newton step on the remaining dofs as outlined in \cref{alg:active_set}.
	The non-convexity of $\mathcal{E}$ leads to catastrophic convergence behavior of the active-set algorithm.
	To overcome this, we employ an extrapolation scheme as used in \cite{HeWhWi15}.
	In the residual \eqref{eq:residual}, this replaces the critical terms $\varphi^2$ in the degradation function $g$ by a prediction $\tilde{\varphi}^2$, computed via linear extrapolation from the solutions of the previous two time-steps.

	\begin{algorithm}[H]
		Repeat for $k=0, \dots$ until $\mathcal{A}_k = \mathcal{A}_{k-1}$ and $\Vert \widetilde{R}_k \Vert \leq \varepsilon_{as}$:
		\begin{algorithmic}[1]			
			\State Determine active-set $\mathcal{A}_k = \{i \mid (M^{-1} R_k)_i + c \ (U_k - U^{old})_i > 0 \text{ and } i \sim \varphi \}$			
			
			\State Fix $\ddphi_k = 0 \quad \forall i \in \mathcal{A}$		
			
			\State Newton-step $U_{k+1} = U_k + G_k^{-1} R_k$ on $\mathcal{A}_k^c$
			
			\State Check convergence
		\end{algorithmic}
		\caption{Primal-dual active-set}
		\label{alg:active_set}
	\end{algorithm}

	The nonlinear residual at the evaluation point $(u, \varphi) \in V \times W_{in}$ is given by 
	\begin{align}\label{eq:residual}
		\begin{split}
			R(u, \varphi)(\du, \dphi) &=
			\left(  g (\tilde\varphi) \sigma^+(u), \strain(\du) \right) 
			+ \left(\sigma^-(u), \strain(\du) \right)
			+ \left(  \partial_\varphi g(\varphi) E_s^+(\strain(u)), \dphi \right)
			\\
			&+ G_c \left( -\frac{1}{\varepsilon} ( 1-\varphi, \dphi ) + \varepsilon \left( \nabla \varphi, \nabla \dphi \right) \right),
		\end{split}
	\end{align}
	with test-functions $(\du, \dphi) \in V \times W$.
	Its derivative at $(u, \varphi)$ in directions $\ddu$ and $\ddphi$ equals to
	\begin{align}\label{eq:G}
		\begin{split}
			G(u, \varphi)(\du, \dphi)(\ddu, \ddphi) =& 
			\left( g(\tilde\varphi) \partial_u \sigma^+(u)(\ddu) , \strain(\du) \right) 
			+ \left( \partial_u \sigma^-(u)(\ddu) , \strain(\du) \right) 
			+ \left( \partial_\varphi g(\varphi) \partial_u(E_s^+(\strain(u)))(\ddu) , \dphi \right)
			\\
			&
			+ \left( \partial^2_{\varphi\varphi} g(\varphi) E_s^+(\strain(u)) \ \ddphi , \dphi  \right) 
			+ G_c \left( \frac{1}{\varepsilon} (\ddphi , \dphi) 
			+ \varepsilon (\nabla \ddphi , \nabla \dphi)  \right).
		\end{split}
	\end{align}
	The partial derivatives wrt $\varphi$ of the first term in $R$ vanish
	{since} %
	we have replaced $\varphi$ by an extrapolation $\tilde\varphi$.
	The corresponding discrete versions at the $k$-th step of \cref{alg:active_set} are denoted by $R_k$ and $G_k$.

	\subsection{Linear Solvers}
	\label{sec:linearsolver}

	The crucial step in the solution of \cref{eq:euler_lagrange} by the active-set method is the solution of the linear system $G_k \cdot \ddU = R_k$ on $\mathcal{A}_k^c$.
	The system to be solved has the following {block} structure
	\begin{equation*}
		\begin{pmatrix}
			G_{uu} & 0 \\
			G_{\varphi u} & G_{\varphi\varphi}.
		\end{pmatrix}
		\cdot
		\begin{pmatrix}
			\ddu \\ \ddphi
		\end{pmatrix}
		=
		\begin{pmatrix}
			R^u \\ R^\varphi
		\end{pmatrix}
		.
	\end{equation*}

	A class of very powerful methods to solve linear systems arising from discretized PDEs are the 
	multigrid methods.
	These aim to reduce and solve the problem 
	{on}
	a hierarchy of nested spaces.
	Two frequently used types of multigrid solvers are the Algebraic Multigrid (AMG) and the Geometric Multigrid (GMG).
	Within the GMG approach, nested spaces are constructed by a hierarchy of meshes, i.e., it requires information about the geometry.
	On the other hand, AMG solvers work solely on the given sparse-matrix.
	There, hierarchies are extracted from the connectivity graph of the matrix.
	Since AMG only needs the linear system to be solved as an input, they are often considered to be black-box solvers.
	This is, however, a huge oversimplification, as there are many different AMG methods and options tailored to specific classes of PDEs.
	For more details on AMG, GMG and multigrid in general, we refer to \cite{Ha85,Br93,St01,HaLa02,TrOoSc01}.
	
	To solve these equations, we use two different approaches.
	First, we use the matrix-free multigrid framework from deal.II extended to handle nonlinear equations as described in our previous work \cite{JoLaWi19}.
	This uses a monolithic geometric multigrid for the whole system.
	The special structure of the matrix is utilized inside the smoother, i.e., we use
	\begin{equation*}
		P^{-1}_{GMG} \coloneqq GMG(G) \mbox{ with block-diagonal smoother } \mathcal{S} \coloneqq 
		\begin{pmatrix}
			S(G_{uu}) & 0 \\
			0 & S(G_{\varphi\varphi})
		\end{pmatrix}
		,
	\end{equation*}
	with Chebyshev-Jacobi smoothers $S$ acting on the single blocks.
	These smoothers are frequently used in the context of matrix-free methods.
	Contrary to {standard} smoothers like Gauss-Seidel, 
	{Chebyshev-Jacobi smoothers}
	only require matrix-vector multiplications and scaling with the diagonal.
	Both operations can be done in a matrix-free fashion, which we will describe later on.
	Chebyshev-type methods work on a given eigenvalue spectrum $[ a, b ]$.
	Depending on the choice of this interval, we can use the Chebyshev method as smoother but also as a solver.
	A solver should handle the whole spectrum of the operator.
	Hence, we want to have $\lambda_{min} \leq a < b \leq \lambda_{max}$ with the extremal eigenvalues $\lambda_{min}$ and $\lambda_{max}$ of $G_{uu}$ resp. $G_{\varphi\varphi}$.
	This is required on the coarsest grid of the multigrid algorithm, where we want to solve the system more accurately.
	The main intention of the Chebyshev-Jacobi method, however, is smoothing.
	Here, we are not interested in treating the whole spectrum, but only the high-frequency parts.
	A choice of $[a,b] := [c \ \lambda_{max}, C \ \lambda_{max}]$ with $ c \leq 1 $ and $ C \geq 1 $ is frequently used.
	In our applications, we use $c = 0.24$ and $C = 1.2$.
	This involves knowledge about $\lambda_{max}$ (and also $\lambda_{min}$ for the solving part).
	For symmetric positive definite matrices, CG can be used to obtain an estimate of both eigenvalues.
	Within PFF, symmetry of $G_{uu}$ and $G_{\varphi\varphi}$ follows from the second derivative of the energy functional.
	We could not rigorously show positive definiteness for these blocks, but numerical studies did not give us contradicting results.
	
	In a second approach, we compare our MF-GMG preconditioner to the block-diagonal AMG preconditioner presented in \cite{HeWi18}.
	The linear system is then preconditioned by
	\begin{equation*}
		P^{-1}_{AMG} \coloneqq
		\begin{pmatrix}
			AMG(G_{uu}) & 0 \\
			0 & AMG(G_{\varphi\varphi})
		\end{pmatrix}
		,
	\end{equation*}
	with Trilinos/MueLu AMG \cite{HeBaHo05,BeGlHu19} solvers for the diagonal blocks.
	In both cases, we use GMRES as outer solver and use $P^{-1}_{AMG}$ resp. $GMG(G)$ as preconditioners.

	\subsection{Matrix-Free Representation}

	The main idea of matrix-free approaches is to avoid assembling and storing sparse (but huge) matrices.
	This is reasonable 
	{since}
	iterative solvers like Conjugate Gradients (CG) or Generalized Minimal Residual (GMRES) do not require explicit knowledge about the entries of the matrix.
	Rather, these solvers only need the result of the matrix-vector multiplication (MV).
	We combine the assembly and sparse matrix-vector multiplication (SpMV) by
	$$G \cdot \ddU = \sum_{k=1}^{n_e} C^T P_k^T G_k (P_k C \ddU),$$
	with the number of elements $n_e$, possible constraints $C$, element-wise global-to-local mapping $P_k$ and the local stiffness matrices $G_k$ corresponding to \eqref{eq:G}.
	This is already a valid matrix-free matrix-vector product (MFMV).
	However, without further optimizations, its performance cannot be expected to be anywhere near to the classical SpMV performance.
	These optimizations are described thoroughly in \cite{KrKo12}, and are included in the deal.II library.
	The main point is to utilize the transformation of the local stiffness matrices $G_k$ to the unit element and abuse the tensor-product structure there.
	In our previous work \cite{JoLaWi19}, we illustrate this in more detail for the nonlinear terms appearing in the PFF problem.

	\subsection{Derivatives of the Eigensystem}
	\label{subsec:DerivativeEigensystem}
	
	In the case of Miehe-type splittings of the stress tensor, we require the derivatives of $\sigma^+$ and $\sigma^-$ in direction $\ddu$, which depend on the eigensystem of the strain tensor.
	In the following, we illustrate how to compute these quantities for arbitrary (but small) dimensions.
	The following derivation holds for real symmetric matrices, and in particular for the solid strain tensor.
	
	We assume, that the matrix $A$ and $\delta A$, as well as the eigenvector and eigenvalue of $A$, denoted by $\lambda$ and $v$, are available such that
	\begin{equation}\label{eq:eigensystem}
		A v = \lambda v, \quad \mbox{ with } \norm{v} = 1,
	\end{equation}
	holds.
	In our application, the direction $\delta A$ is given in terms of the displacement test functions.
	The eigensystem of $A$ itself can easily be computed explicitly in $2d$ and $3d$, e.g., with the help of some computer algebra system.
	For higher dimensions, more sophisticated methods are available.
	Taking the derivative of \eqref{eq:eigensystem}, we get by the product rule
	\begin{equation}\label{eq:lin_eigensystem}
		\delta A \cdot v + A \cdot \delta v = \delta \lambda \ v + \lambda \ \delta v.
	\end{equation}
	Forming the $\ell_2$-scalar-product of \eqref{eq:lin_eigensystem} with $u$ yields
	\begin{equation*}
		\sprod{\delta A \cdot v}{v} + \sprod{A \cdot \delta v}{v} = \sprod{\delta \lambda \ v}{v} + \sprod{\lambda \ \delta v}{v}.
	\end{equation*}
	Using the symmetry of $A$ and the normalization $\sprod{v}{v} = 1$, 
	{we can reduce it to}
	\begin{equation}\label{eq:lin_lambda}
		\sprod{\delta A \cdot v}{v} = \delta \lambda.
	\end{equation}
	{Since} 
	$\delta A$ and $v$ are given, we can now compute the derivative of $\lambda$ by the formula above.

	Computing the derivative of the eigenvector $v$ is slightly more involved.
	We start by collecting the terms $\delta v$ in \eqref{eq:lin_eigensystem}:
	\begin{equation}\label{eq:eigensystem:L}
		(A - \lambda I) \cdot \delta v = \delta \lambda \ v - \delta A \cdot v.
	\end{equation}
	Unfortunately, we cannot simply invert $(A - \lambda I) \eqqcolon L$ to obtain $\delta v$
	{since}
	$\lambda$ is an eigenvalue of $A$. Hence, $L$ is singular by definition.
	This is seen easily from \eqref{eq:eigensystem}, which can be rearranged to $L \cdot v = 0$.
	Therefore, we need the concept of a pseudo-inverse (or left-inverse), see e.g. \cite{BeGr80}.
	The pseudo-inverse $L^\dagger$ has the properties (among others) that $L^\dagger L = I$ and $L^\dagger w = 0 \quad \forall w \in \mathcal{N}(L) = \mbox{span}\{ v \}$.
	Utilizing this, we can apply $L^\dagger$ to \eqref{eq:eigensystem:L}, and get
	\begin{equation*}
		L^\dagger L \cdot \delta v = L^\dagger \delta \lambda \ v - L^\dagger \delta A \cdot v,
	\end{equation*}
	which finally reduces to
	\begin{equation}\label{eq:lin_v}
		\delta v = - L^\dagger \delta A \cdot v.
	\end{equation}	
	With \eqref{eq:lin_lambda} and \eqref{eq:lin_v} we can now compute the linearized eigensystem of $A$.

	If the eigensystem $(\lambda_i, v_i)$ of $A$ is known, we can decompose $A = V D V^T$ with the matrix of eigenvectors $V = (v_1, \dots, v_d)$ and diagonal matrix $D = \mbox{diag}(\lambda_1, \dots, \lambda_d)$.
	Furthermore, we have $L = V (D - \lambda I) V^T$.
	Then, the pseudo-inverse can be computed as
	\begin{equation*}
		L^\dagger = V D^\dagger V^T,
	\end{equation*} 
	with $D^\dagger \coloneqq 
	\begin{cases} 
		0 & \mbox{if } \lambda_i - \lambda = 0 \\
		\frac{1}{\lambda_i - \lambda} & \mbox{if } \lambda_i - \lambda \neq 0. \\
	\end{cases}$
	
	This machinery can be applied to compute $\strain^+$ and its derivative for arbitrary dimensions, see \cref{alg:strain_plus} for the detailed steps.
	The derivative of $\sigma^+(e)$ follows by
	$$\partial_u \sigma^+(u)(du) = \lambda \begin{rcases} \begin{dcases} 0 & \mbox{if } \mbox{tr}(e(u)) < 0 \\ \mbox{tr}(e(du)) & \mbox{else} \end{dcases} \end{rcases} I + 2 \mu \partial_e (e^+) (e(du)).$$
	Furthermore, each of the steps can be implemented using SIMD parallelism.
	The (hidden) "if" statement in line (10) can be handled via masking as explained in \cref{sec:simd}.

	\begin{algorithm}
		\begin{algorithmic}[1]
			\State \textbf{Input:} $A, \delta A \in \mathbb{R}^{d \times d}, d = 2,3$
			\State \textbf{Output:} $\partial_A A^+(A)(\delta A) \in \mathbb{R}^{d \times d}$
			\Comment{Derivative of $A^+$ at $A$ in direction $\delta A$}
			\State
			\State Decompose $A = V D V^T$
			\State
			\For{$i \in \{1, \dots,  d \}$}
			\State $L = A - \lambda_i I$
			\State $D^\dagger \coloneqq \mbox{diag}(\frac{1}{\lambda_1 - \lambda_i}, \dots, \frac{1}{\lambda_d - \lambda_i})$
			\State $L^\dagger = V D^\dagger V^T$
			\State $dv_i = - L^\dagger \delta A \cdot v_i$
			\State $d\lambda_i = v_i^T \ \delta A \ v_i$
			\State $d\lambda_i^+ = \begin{cases}
			0  & \mbox{if } \lambda_i < 0 \\
			d\lambda_i & \mbox{if } \lambda_i \geq 0
			\end{cases}$
			\EndFor
			\State
			\State $\delta V = ( dv_1, \dots, dv_d )$
			\State $\delta D^+ = \mbox{diag}(\delta \lambda^+_1, \dots, \delta \lambda^+_d)$
			\State
			\State $A^+ \coloneqq V D^+ V^T$
			\State
			\State \Return $\partial A^+(A)(\delta A) \coloneqq \delta V D^+ V^T + V \delta D^+ V^T + V D^+ \delta V^T$
		\end{algorithmic}
		
		\caption{
		}
		\label{alg:strain_plus}
	\end{algorithm}

\section{Parallel Implementation}
\label{sec:parallel}

	The finite element library deal.II as well as our implementation are parallelized on multiple levels.
	First of all, the problem is split across multiple CPUs such that each of them only stores and handles a small subset of the whole problem (distributed parallelization).
	Necessary communication between different CPUs is handled via the message passing interface (MPI).
	Furthermore, each CPU runs explicitly vectorized code, that is, multiple data elements are handled with the same CPU instruction (SIMD).
	The next subsections explain these concepts in more detail.

	\subsection{Distributed Parallelization}
	\label{sec:distributed}
	
	When dealing with large problems, it may no longer be possible for a single CPU to store all the required infrastructure in memory.
	Hence, the computational problem is split into smaller chunks which are distributed among multiple CPUs.
	In the context of FEM, this is typically done by assigning each CPU some parts of the mesh.
	Such a partitioning can be achieved by using libraries like p4est \cite{BuWiGh11}.
	Each CPU then handles all the dofs associated to its part of the mesh (owned dofs).
	Furthermore, it may occasionally require knowledge about additional dofs which are owned by another CPU (ghost dofs), e.g. if we want to evaluate solutions close to the interface between different CPUs.
	These information needs to be shared among multiple CPUs, which is done via the message passing interface (MPI).
	This provides methods for communication in distributed networks.
	deal.II provides an easy to use infrastructure for these kind of tasks; see, e.g., \cite{BaBuHe11}.

	\subsection{Vectorization}
	\label{sec:simd}
	
	Within the matrix-free framework, deal.II makes use of the CPUs vector instruction set.
	These operations perform single instructions on multiple data (SIMD).	
	Consider the simple implementation of adding two vectors: $r[i] = a[i] + b[i] \quad i = 0, \dots 3$.
	This would take the CPU to issue one \textit{add} instruction per vector entry, i.e., a total of $4$ instructions (ignoring CPU instructions required for the loop, assignment, etc.).
	However, modern CPUs are capable of performing the same operations on multiple elements at once (not to be confused with multithreading~!).
	Hence, only $1$ \textit{add} instruction is required for the simple example above, resulting in a hypothetical speedup of $4$.
	Depending on the CPU, different vectorization sizes are possible.
	The most recent CPUs support $128$ bit, $256$ bit and $512$ bit instructions.
	Hence, $2, 4$ and $8$ double precision values or $4, 8$ and $16$ single precision values may be treated simultaneously.
	
	There are several important points to consider here.
	First of all, vectorization requires the same work to be done on all entries.
	Hence, code with branching, e.g., \textit{if} statements with a condition depending on the values cannot be efficiently vectorized.
	Workarounds usually involve executing both branches for each vector element and merging them via bit masks, see \cref{alg:if_then_else}.
	These type of instructions are called "blend".
	This is mainly used for simple \textit{if-then-else} assignments, i.e., $x[i] = \begin{cases} a[i] \quad &\text{if } a[i] > b[i] \\ 0 &\text{else} \end{cases}$, where the additional overhead of executing all branches and combining the results does not dominate the benefits of using vectorized code.

	\begin{algorithm}[H]
		\begin{algorithmic}[1]
			\State Input: $a, b \in \mathbb{R}^v$, $v = 1, 2, 4,$ or $8$
			\State Output: $x[i] = \begin{cases} a[i] \quad &\text{if } a[i] > b[i] \\ 0 &\text{else} \end{cases}$
			\State $mask = x < y$ 
			\Comment{a bitmask, with all bits of element $i$ set to $1$ if $x[i] < y[i]$, and or $0$ otherwise}
			\State $a0 = a$ \& mask  \Comment{bitwise 'and' operation. Each entry equals to either the entry of $a$ or $0$}
			\State $b0 = b$ \& !mask \Comment{similar to above: each entry is either part of b (condition was false) or $0$}
			\State $result = a0 \ | \ b0$ \Comment{combine the two branches by a bitwise 'or' operation}
		\end{algorithmic}
		\caption{Vectorized if-then-else statement using bit masks.}
		\label{alg:if_then_else}
	\end{algorithm}
	
	An additional drawback of vectorization is the difficulty of implementation and portability.
	For simple algorithms, up-to-date compilers can automatically vectorize the loop.
	Unfortunately, this auto-vectorization ceases to work for more complex scenarios.
	This leaves the developer to explicitly call SIMD instructions of the form \textit{\_mm512\_add\_pd}, which makes the implementation more tedious.
	Usually, wrapper libraries are used to improve code readability and portability.
	{For instance,}
	the VectorizedArray class in deal.II or Vc \cite{KrLi12} among others.
	Recently, there has also been effort to include vectorization capabilities into the C++ standard.
	
	Within the deal.II matrix-free implementation, vectorization is done over elements.
	That is, each CPU instruction acts on multiple elements at once.
	Several implementational details are described in \cite{KrKo12,KrKo19,KrWa18}.
	
	We now investigate the influence of different vectorization levels and compiler settings on the performance of the MFMV.
	All of these test are carried out on a single core, using the single-edge notched shear test described in \cref{sec:shear}.
	We used $Q_1$ elements on a $6$-times refined grid, but the observed speedup is valid also for different settings, as long as the problem size is not too small.
	Compilation was done using GCC-9.2.0 using the O3 optimization level.
	The test ran on an Intel Haswell Xeon E5-2630v3 with $2.4$ GHz clock speed, which did not support the most recent instruction set AVX-512.
	All computations were carried out using double precision floating point numbers.
	The best run out of $100$ is reported.

	Our findings are summarized in \cref{tab:simd}.
	There, "ftree-vectorize" and "fno-tree-vectorize" are the GCC flags to enable resp. disable automatic vectorization.
	$n$-bits denotes explicit vectorization using the given number of bits, i.e. $128$ for SSE and $256$ for the AVX instruction set.
	
	We notice that auto-vectorization does only yield a speedup of $13\%$.
	Explicit vectorization over elements doubles the performance, hence, yielding the best possible speedup using the $128$ bit instruction set.
	Additionally enabling auto-vectorization still increases the performance a little bit.
	This is mainly due to vectorization of simple parts like adding vectors, that are not explicitly vectorized but relying on auto-vectorization.
	Vectorization capabilities of compilers have improved throughout the years.
	However, efficient vectorization requires knowledge about the structure of the computation to be carried out that is beyond the compilers ability to identify, e.g., realizing that the work on most elements is identical.
	
	Using $256$ bit SIMD instructions only improved the results by another $27\%$.
	Possible explanations are the automatic throttling of the CPU speed \cite{In19} or memory bandwidth restrictions, but the actual reasons are unclear.
	Surprisingly, auto-vectorization interferes negatively with explicit vectorization in any case.
	
	\begin{table}
		\centering	
		\begin{filecontents*}{tables/slit_simd_gcc.csv}
idx_folder run time type idx_pattern bit_width name args speedup 
0 86 0.030400989 double 0 64 o3 {$O3$, fno-tree-vectorize} 1.0 
1 36 0.026873948 double 0 64 auto {$O3$, ftree-vectorize} 1.1312438723182763 
3 42 0.014892786 double 0 128 expl_auto {128 bits, ftree-vectorize} 2.041323161428627 
2 63 0.014447336 double 0 128 explicit {128 bits, fno-tree-vectorize} 2.1042626128443334 
3 22 0.013852043 double 0 256 expl_auto {256 bits, ftree-vectorize} 2.194693519215902 
2 76 0.013409956 double 0 256 explicit {256 bits, fno-tree-vectorize} 2.2670461409418494 
\end{filecontents*}
\pgfplotstabletypeset[
			columns={args, speedup},
			columns/args/.style={column name=Arguments, string type, column type = {l}},
			columns/speedup/.style={column name=Speed-Up, precision = 2, zerofill},
			every head row/.style={before row=\toprule,after row=\midrule},
			every last row/.style={after row=\bottomrule},
		]{tables/slit_simd_gcc.csv}

		\caption{
			Effect of vectorization within the matrix-free operator evaluation.
			Baseline is given by a standard release build using the $O3$ optimization flag without vectorization.
		}
		\label{tab:simd}
	\end{table}

\section{Numerical Examples}
\label{sec:numericresults}

	In the upcoming sections, we present our numerical examples.
	We focus on a pair of well-known test cases for fracture propagation: the single-edge-notched shear and tension test.
	These have been well studied in the literature, such that we can compare our findings to those reported by other groups.
	Results regarding the performance of our solvers and their behavior under $p$-refinement follows in \cref{sec:compare}.

	\subsection{Example 1: Single-Edge Notched Tension Test}

	We start with the single-edge notched tension test.
	This test case has been studied by other groups as well, 
        see e.g. \cite{MiWeHo10,MiHoWe10,AmGeDe15,Bo12,HeWhWi15}.
	
	\subsubsection{Description}

	The domain of interest is a unit square with a thin slit as shown in \cref{fig:tension:geometry}.
	Due to the discontinuity in the geometry, we have two sets of dofs along the slit at the same locations.
	This is indicated by the double dot at $(0,0.5)$.
	Fracture propagation is driven by an increasing displacement pulling at the top boundary, 
	i.e., $u = (0, t \cdot du) \mbox{ on } \Gamma_{top}$.
	The specimen is fixed at bottom by imposing $u = (0,0) \mbox{ on } \Gamma_{bottom}$.

	\begin{figure}[H]
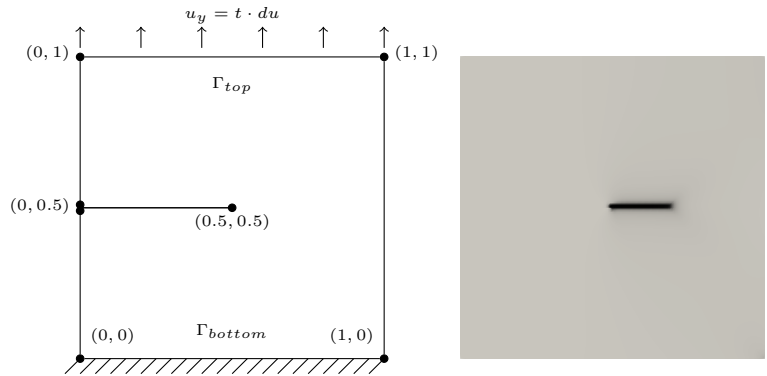

		\centering
		\include{img/slit_tension}		
		\caption{
			(\textbf{a}) Geometry for the single-edge notched tension tests (values given in $mm$.).
			A predefined slit is imposed in the geometry of the domain.
			(\textbf{b}) The computed fracture moves on a straight line towards the right boundary.
		}
		\label{fig:tension:geometry}
	\end{figure}

	\begin{Remark}
		Due to symmetry, it would be sufficient to fix the $y$-component of $u$ at $\Gamma_{bottom}$ and $\Gamma_{top}$ for the tension test.
		This yields the same results as using the boundary conditions specified above.
	\end{Remark}

	The material parameters for the tension test are shown in \cref{tab:tension:material}.
	All computations were done using $Q_1$-elements for displacement and phase-field on uniform quadrilateral grids.
	\Cref{tab:tension:dofs} lists the number of elements and dofs for each level of refinement $\ell$. 

	\begin{table}[H]
		\renewcommand{\arraystretch}{1.3}
		\centering
		\begin{tabular}{l|cccccc}
			Parameter & $\lambda$ & $\mu$     & $G_c$               & $\kappa$   & $dt$ & $du$      \\\hline
			Value     & $121.15$  & $80.77$   & $2.7 \cdot 10^{-3}$ & $10^{-10}$ & $1$  & $10^{-5}$ \\\hline
			Unit      & $kN/mm^2$ & $kN/mm^2$ & $kN/mm$             &            & $s$  & $mm$     
		\end{tabular}
		\caption{Material parameters and configuration for the single-edge notched shear/tension test cases.}
		\label{tab:tension:material}
	\end{table}

	\begin{table}[H]
		\renewcommand{\arraystretch}{1.3}
		\centering
		\summary{tables/slit_tension_none.csv}
		\caption{Values and parameters for the different refinement levels of the single-edge notched shear/tension tests.}
		\label{tab:tension:dofs}
	\end{table}

	\subsubsection{Discussion}

	First, we perform a refinement study for the single-edge notched tension test.
	To this end, we fix all parameters and perform global mesh refinement.
	In particular, the phase-field parameter is fixed to $\varepsilon = 4 \cdot 10^{-3}$ as used in \cite{AmGeDe15}.

	The resulting load-displacement curves are shown in \cref{fig:tension:refinement} for the anisotropic case (Miehe splitting).
	Whereas the solutions on the coarse meshes deviate a lot, the solutions on the fine meshes are almost indistinguishable showing the convergence of our implementation.
	In particular, the deviation to the solution on the finest mesh seems to decrease by approximately $0.5$ in every refinement step.
	The splitting does not seem to have much impact on this simulation. 
	{Indeed,} the isotropic results agree almost perfectly with the anisotropic findings.

	In comparison to the results presented in \cite{AmGeDe15}, the resulting maximum load is slightly higher in our computations.
	Furthermore, the decay of the loading force seems to take longer, i.e. the amount of displacement when the load hits zero is roughly $7\cdot 10^{-3}$ in our case, whereas it is slightly less than $6 \cdot 10^{-3}$ in \cite{AmGeDe15}.
	However, we would like to point out that our computations used significantly more elements.

	\begin{figure}[H]
		\centering
		\includegraphics{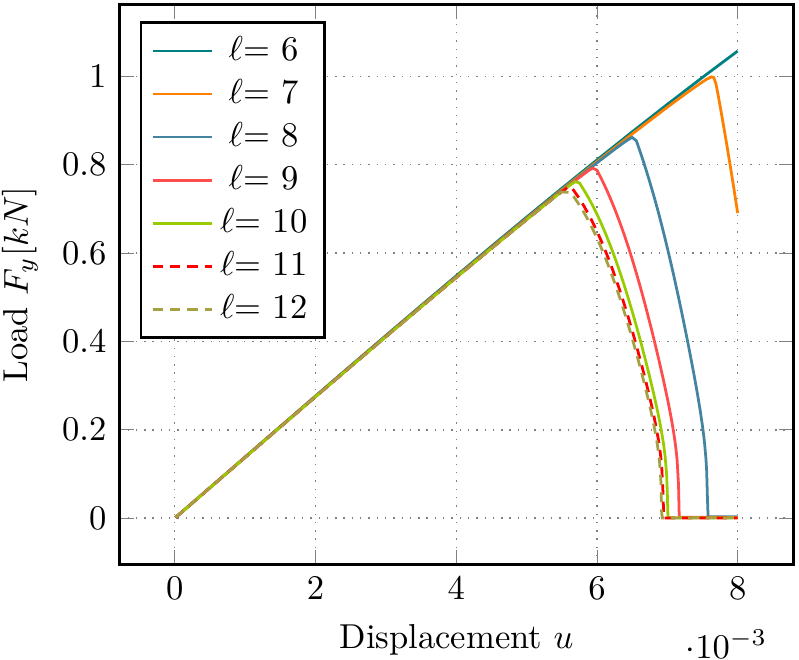}
		\caption{
			Refinement studies for the single-edge notched tension test using Miehe splitting.
			All other parameters are fixed.
			In particular, $\varepsilon$ is set to $\varepsilon = 4 \cdot 10^{-3}$.
		}
		\label{fig:tension:refinement}
	\end{figure}

	\subsection{Example 2: Single-Edge Notched Shear Test}
	\label{sec:shear}

	The following test scenario is closely related to the previous tension test.
	Both have the same geometry and material parameters, but the boundary conditions are flipped.

	\subsubsection{Description}

	The computational domain is again given by a unit square, with the same slit imposed as in the previous test case.
	In this scenario, however, the applied displacement conditions moves {from} the top boundary to the right, i.e., $u = (t \cdot du) \mbox{ on } \Gamma_{top}$.
	Again, the 
	{body remains}
	fixed at the bottom by $u = (0,0) \mbox{ on } \Gamma_{bottom}$.

	\begin{figure}[H]
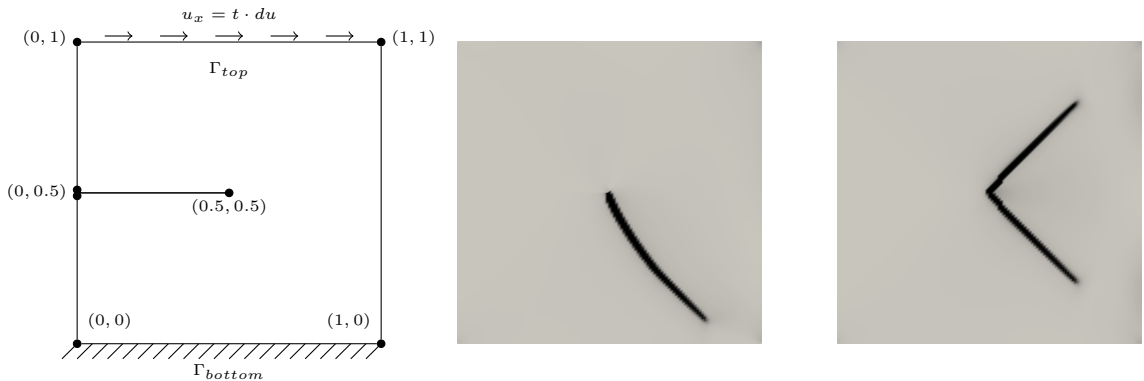

		\centering
		\include{img/slit_shear}
		\caption{
			(\textbf{a}) Geometry for the single-edge notched shear tests (values given in $mm$).
			(\textbf{b}) If we apply Miehe-splitting, the fracture moves along a curved path towards the bottom of the 
			{body}
			.
			(\textbf{c}) In the isotropic case (no splitting), two symmetric crack branches emerge 
			.
		}
		\label{fig:shear:geometry}
	\end{figure}

	\subsubsection{Discussion}

	In the isotropic case of the single-edge notched shear test in \cref{fig:shear:loaddisplacement}, the load-displacement curves do not coincide with the observations in \cite{AmGeDe15} but look more similar to their results with Amor-splitting. 
	However, the observed crack pattern is similar, i.e., two symmetric branches appear as seen in \cref{fig:shear:geometry} (right).
	We like to point out that neither the Miehe-splitting nor the isotropic model lead to physically sound results, see e.g. \cite{AmGeDe15} and 
also the early work \cite{MiWeHo10}.
	Nonetheless, one clearly observes convergence as in the tension test.
	
	Contrary to the tension test before, the resulting load-displacement curves show lots of small oscillations.
	This is partly also observed in \cite{AmGeDe15}, although the oscillations there are seemingly smaller.
	The results with Miehe splitting seem match better, but the maximum load and displacement are lower than those reported by \cite{AmGeDe15}.
	
	\Cref{fig:shear:loaddisplacement:ho} shows the test results for different polynomial degrees $p$.
	The refinement level corresponds to $\ell = 7$, all other parameters are the same as before.
	We immediately notice, that the small oscillations are gone for the high-order simulations.
	Furthermore, we also see convergence with respect to $p$.
	The number of dofs for degree $p$ can be obtained from \cref{tab:slit:dofs:p}.

	\begin{table}[H]
		\renewcommand{\arraystretch}{1.3}
		\centering
		\begin{filecontents*}{tables/slit_shear_p.csv}
lvl dofs eps h run cells degree 
7 50115 0.004 0.0110485435 0 16384 1 
7 198531 0.004 0.0110485435 0 16384 2 
7 445251 0.004 0.0110485435 0 16384 3 
7 790275 0.004 0.0110485435 0 16384 4 
\end{filecontents*}
\pgfplotstabletypeset[
		columns={lvl, degree, dofs},
		columns/lvl/.style={column name=$\ell$},
		columns/degree/.style={column name=$p$},
		columns/dofs/.style={sci, zerofill, precision = 1},
		every head row/.style={before row=\toprule,after row=\midrule},
		every last row/.style={after row=\bottomrule},
		]{tables/slit_shear_p.csv}
		\hspace{1cm}
		\begin{filecontents*}{tables/slit_scaling_p.csv}
lvl dofs eps h run cells degree 
9 790275 0.004 0.00276213586 0 262144 1 
9 3153411 0.004 0.00276213586 1 262144 2 
9 7089411 0.004 0.00276213586 2 262144 3 
9 12598275 0.004 0.00276213586 3 262144 4 
\end{filecontents*}
\pgfplotstabletypeset[
		columns={lvl, degree, dofs},
		columns/lvl/.style={column name=$\ell$},
		columns/degree/.style={column name=$p$},
		columns/dofs/.style={sci, zerofill, precision = 1},
		every head row/.style={before row=\toprule,after row=\midrule},
		every last row/.style={after row=\bottomrule},
		]{tables/slit_scaling_p.csv}
		
		\caption{Dof counts for varying $p$ for the single-edge notched shear/tension tests.
		\textbf{(a)} mesh-refinement level $\ell = 7$, \textbf{(b)} $\ell = 9$.}
		\label{tab:slit:dofs:p}
	\end{table}
	
	\begin{figure}[H]
		\centering
		\includegraphics{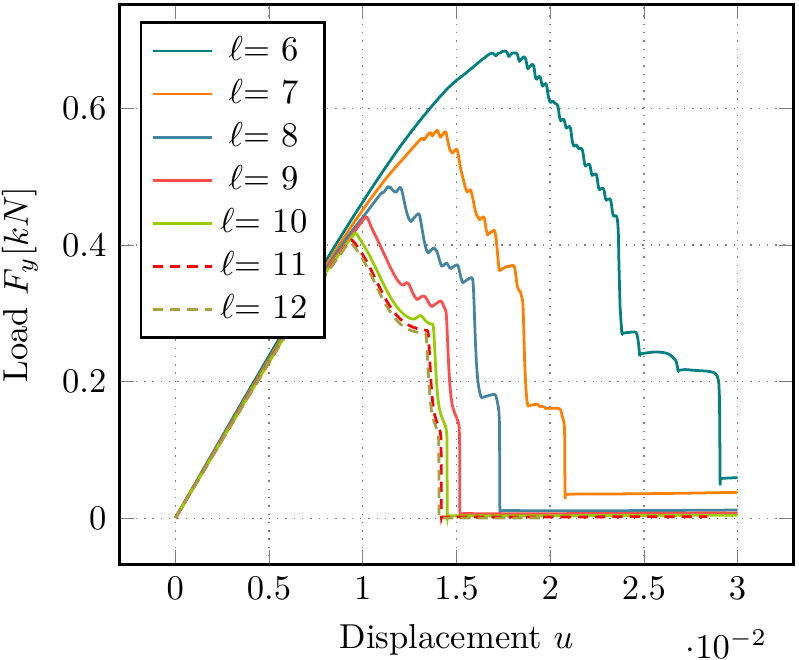}%
		\includegraphics{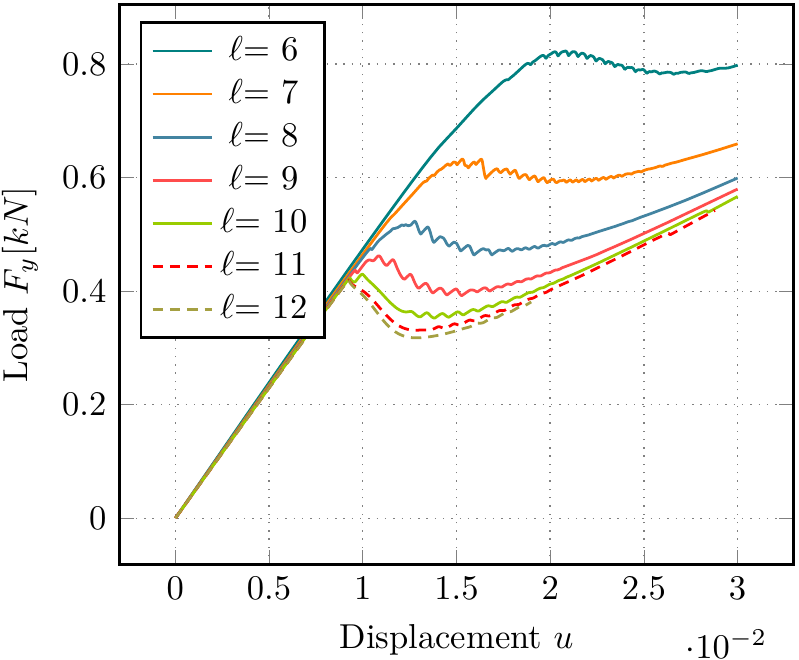}
		
		\caption{
		Refinement studies for the single-edge notched shear test.
		(\textbf{a}) Isotropic case (no splitting). (\textbf{b}) Miehe splitting.
		All other parameters are fixed.
		In particular, $\varepsilon$ is set to $\varepsilon = 4 \cdot 10^{-3}$.
		}
		\label{fig:shear:loaddisplacement}
	\end{figure}

	\begin{figure}[H]
		\centering
		\includegraphics{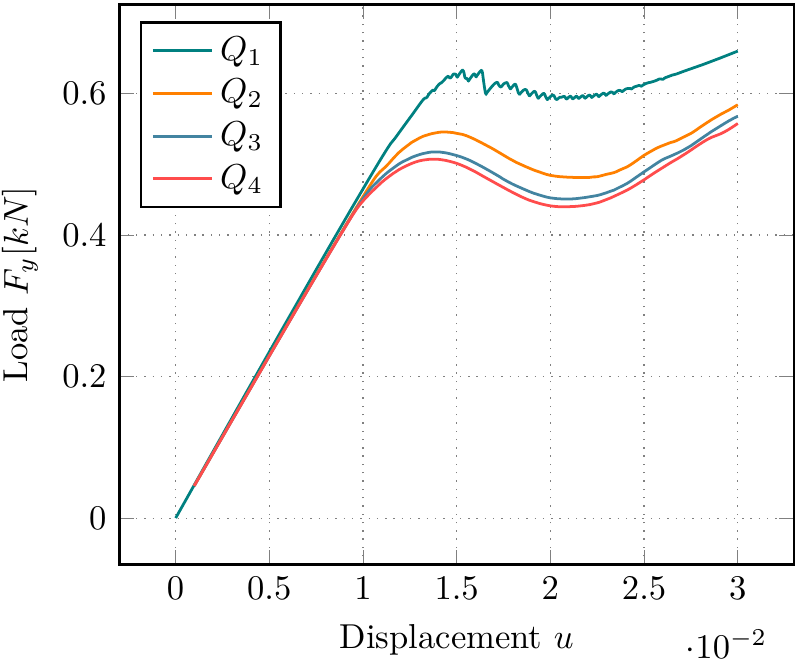}
		
		\caption{
			Different polynomial degrees for the anisotropic single-edge notched shear test.
			Refinement level $\ell$ is fixed at $\ell = 7$ and $\varepsilon = 4 \cdot 10^{-3}$.
		}
		\label{fig:shear:loaddisplacement:ho}
	\end{figure}

\section{Performance Studies}
\label{sec:compare}

	In the upcoming sections, we compare our matrix-free version with the AMG based preconditioner described in \cite{HeWi18} in terms of iteration counts, memory requirements, computational time and scalability.
	We would like to point out that most of the code is shared between the two implementations, i.e., only the parts related to handle the linear system differ.
	This allows for a fair comparison of the results.

	\begin{figure}[H]
		\centering
		\includegraphics{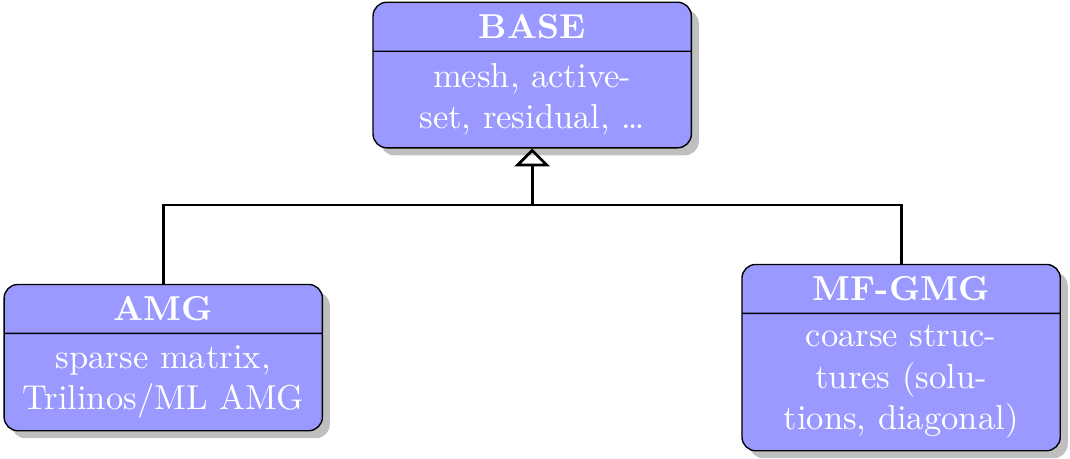}
		
		\caption{Shared code base for the AMG and MF-GMG implementation.}
		\label{fig:code:inheritance}
	\end{figure}

	\subsection{Memory Requirements}

	We start by investigating the memory requirements of the two approaches, which should be an obvious advantage of the matrix-free framework.
	Due to the required multigrid dofs for the geometric multigrid, the dof structures require slightly more memory compared to the AMG version.
	This is more than compensated by larger storage requirements of the sparse matrix and AMG hierarchies compared to the multilevel matrix-free structures in the GMG implementation.
	In particular for high-order polynomial degrees $p$, the sparse matrix gets increasingly dense, leading to huge memory costs.	
	Unlike that, the cost for a dof within the matrix-free approach is almost independent of $p$.

	All these effects are shown in \cref{fig:memory} for varying polynomial degree $p$.
	The line "Matrix (AMG)" refers to the storage cost of the sparse-matrix.
	The entry "Matrix (GMG)" denotes memory requirements for all necessary matrix-free operators, and additionally required data, including those on the coarser levels needed for the geometric multigrid.
	The additional cost of handling the dofs on coarser levels as well is reflected in the difference between "DoFs (AMG)" and "DoFs (GMG)".
	For a degree of $4$, we save approximately a factor of $20$ in terms of memory.
	These results are shown for a $2d$ test case.
	{In 3d,}
	this gap would 
	{even be}
	more prominent.
	Note that you must not get confused by the decrease of the "Matrix (GMG)" line, as we plot the required bytes per dof and not the total memory needed.

	\begin{Remark}
		Data structures of the deal.II library 
		{(matrix-free, GMG, dofs, etc.)}
		are able to report their memory consumption very precisely.
		Unfortunately, this is not the case for the AMG structures of Trilinos/ML.
		Measuring the memory requirements for the AMG using \path{/proc/self/status} did not yield reliable results.
		Hence, we did not include the storage of AMG structures in \cref{fig:memory}.
		Although this makes the comparison slightly unfair, it would only make the results even more favorable for the matrix-free implementation.
	\end{Remark}

	\begin{figure}[H]
		\centering
		\includegraphics{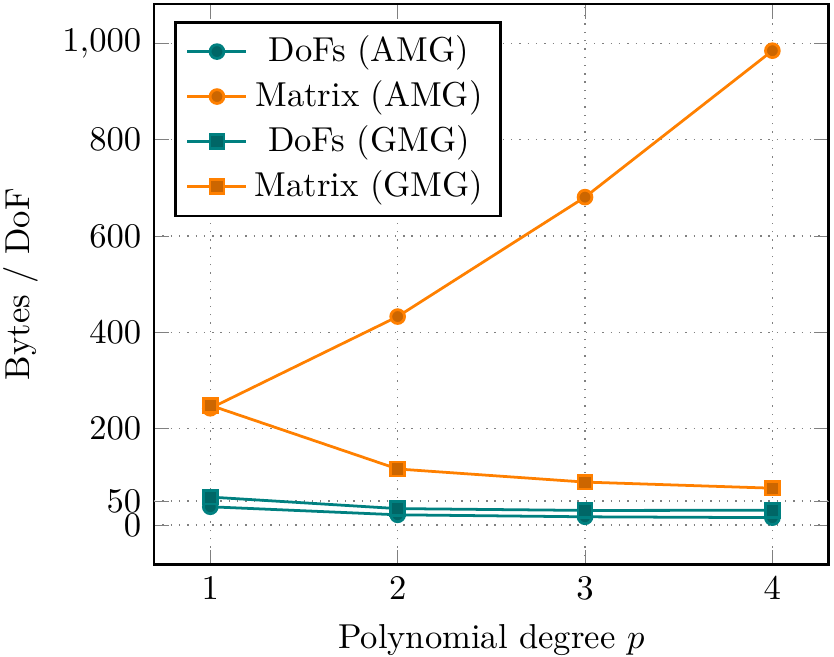}	
		\caption{
				Storage requirements in bytes per dof for varying polynomial degrees.
		}
		\label{fig:memory}
	\end{figure}

	\subsection{Iteration Counts}

	Next, we have a look 
	{at}
	the number of 
	iterations required to solve the arising {linear} systems of equations.
	We would like to point out that these counts are highly dependent on the actual settings used, i.e., which simulation is run, selected number of smoothing steps, coarse level and further parameters.
	For the AMG solver, we stuck to the default configuration by deal.II.
	The GMG uses $\ell = 2$ as coarse grid, i.e. $16$ elements.
	We always use $5$ Chebyshev-Jacobi smoothing iteration sweeps.
	Increasing the number of smoothing steps would decrease the number of iterations.
	Finding a good balance between smoothing quality and iteration count to obtain the best performance is a very challenging task, 
        which we did not consider {in detail} by now.
    In all our simulations, we use relative convergence criteria of $10^{-6}$ for the linear solver and $10^{-8}$ for the active-set strategy.

	A representative comparison of the iteration behavior is given in \cref{fig:compare:iter}, 
	with AMG results 
	{on}
	the left and GMG results on the right side.
	Again, we vary the polynomial degree for the shear test, as this scenario shows more interesting behavior.
	Both approaches show an expected growth in the number of iterations for increasing $p$.
	We also observe that both solvers require more work as the fracture grows.	
	Degree $p=4$ behaves somewhat surprisingly.
	{It}
	requires less steps in the middle of the simulation for {both} the AMG and GMG strategy.
	{The} iteration count jumps up towards the end of the simulation, when the domain is close to total failure.

	\begin{figure}[H]
		\centering
		 \includegraphics{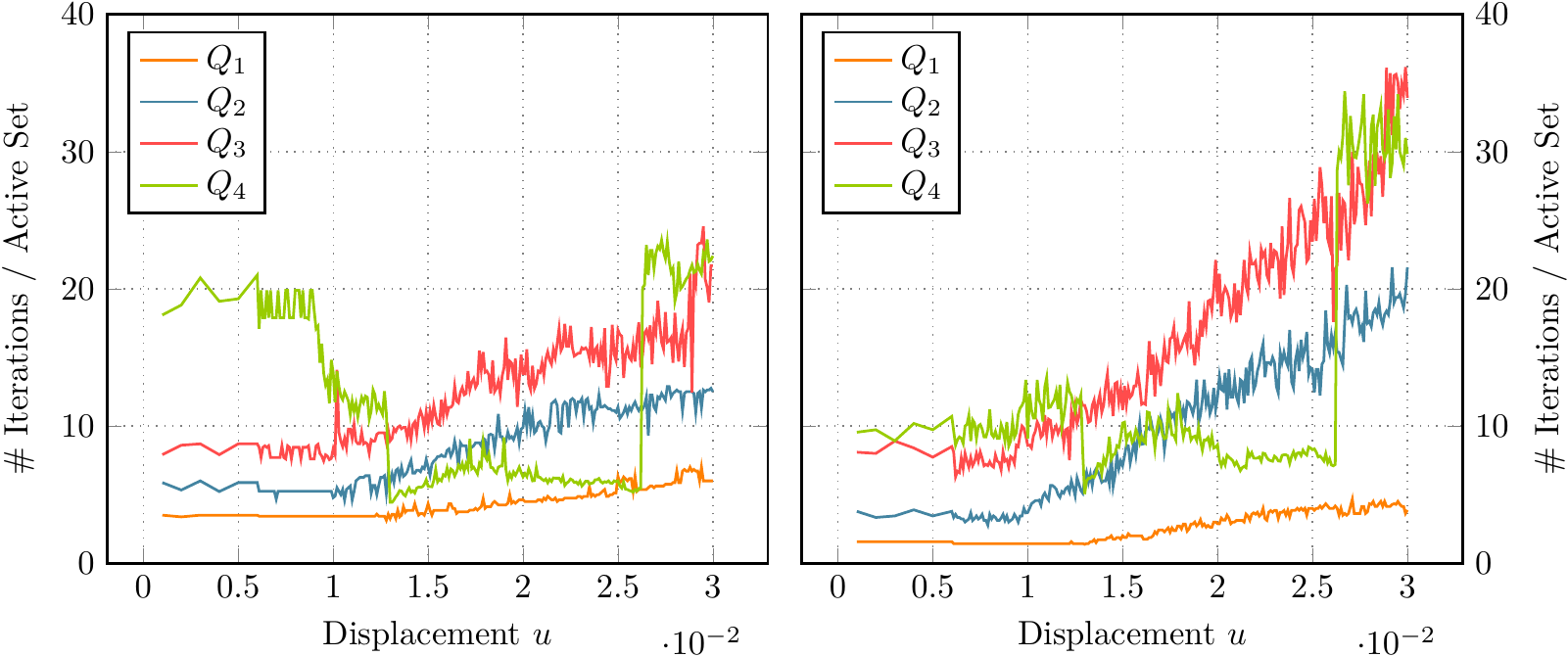}
		
		\caption{
			Iteration counts using the single-edge notched shear test.
			(\textbf{a}) Block-diagonal AMG.  (\textbf{b}) MF-GMG preconditioner.
			Both strategies require more iterations once the domain is close to total failure.
			Mesh refinement corresponds to $\ell = 7$, the length-scale parameter is fixed at $\varepsilon = 4 \cdot 10^{-3}$.
		}
		\label{fig:compare:iter}
	\end{figure}

	\subsection{Performance Analysis}
	
	By avoiding storing the matrix explicitly, we are forced to compute (almost) everything on the fly.
	This makes matrix-vector operations more costly compared to the standard case, where the matrix is at our disposal.
	On the other hand, matrix-free implementations start to outperform standard methods for higher-order elements, in particular in $3d$ computations.

	We continue to investigate the computational performance of the matrix-free approach.
	To this end, we start by comparing the sparse-matrix vector multiplication (SpMV) time with the corresponding matrix-free evaluation (MFMV).
	For all 
	tests, we use explicit vectorization 
	{based on}
	AVX-256 in the matrix-free parts as described in \cref{sec:simd}.
	The most recent instruction set AVX-512 is not supported by any of our CPUs.
	We plot the best run out of $100$.
	
	In \cite{JoLaWi19}, we compared the behavior with respect to $h$-refinement.
	Here, we want to focus on the dependence on the polynomial degree $p$.
	In \cref{fig:compare:vmult}, we observe that SpMV is faster for linear and quadratic elements.
	However, this 
	{only considers} the raw matrix-vector multiplication~!
	In particular, it does not include the time 
	{that is required}
	to assemble the matrix.
	We also notice, that starting from degree $3$, MFMV starts to outperform SpMV even for single evaluations.

	\begin{figure}[H]
		\centering
		\includegraphics{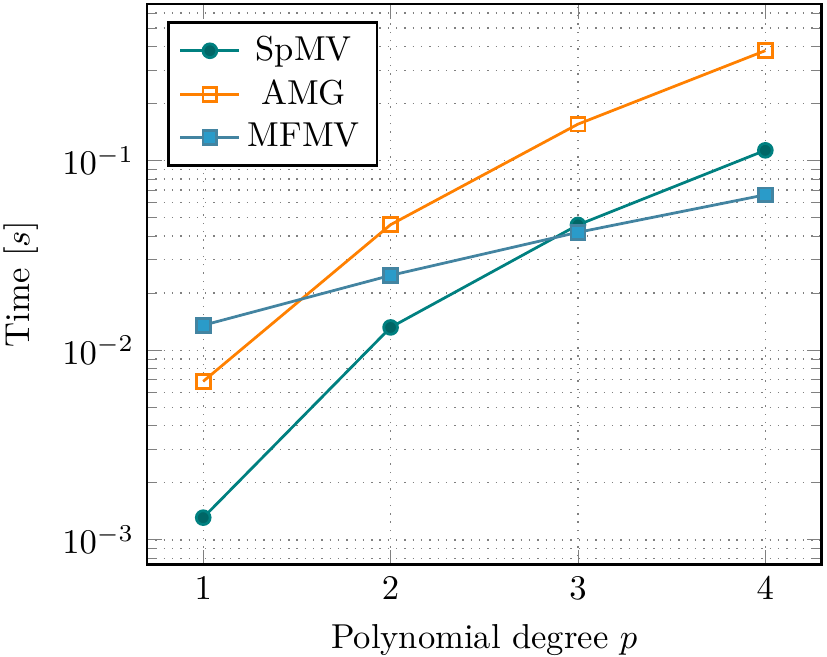}
		\caption{
			Comparison of SpMV and MFMV for varying polynomial degrees.
			Again, $\ell = 7$ and $\varepsilon = 4 \cdot 10^{-3}$.
		}
		\label{fig:compare:vmult}
	\end{figure}
	
	{Let us now look at the perfomance of a full simulation.}
	This is visualized in \cref{fig:compare:time}, where the total time spent in the linear solvers, i.e., GMRES with either AMG or GMG as preconditioners, is plotted.
	Here, only a single simulation of the single-edged notched shear test is performed.
	Possible fluctuations in the timings should be evened out over the length of the simulation.
	
	The observed behavior is very similar to the vmult-times presented before.
	The faster SpMV for linear and quadratic elements is able to overcome the time of assembling the matrix and initializing the AMG solver, thus ending up faster than the matrix-free solver.
	The situation shifts in favor of the MF-GMG preconditioner for polynomial degrees $3$ and higher.

	\begin{figure}[H]
		\centering
		\includegraphics{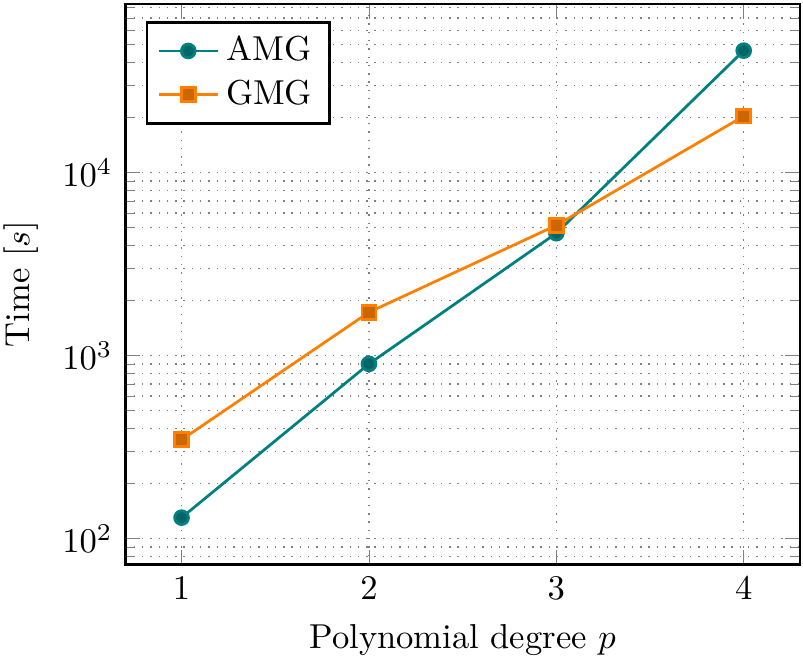}
		\caption{
			Total time spent solving the linear system for a full simulation of the single-edge notched shear test using $16$ cores.
			This includes setup times for the sparse-matrix and AMG solver as well as the eigenvalue computation and setup routines required for the MF-GMG.
			Again, $\ell = 7$ and $\varepsilon = 4 \cdot 10^{-3}$.
		}
		\label{fig:compare:time}
	\end{figure}

	\subsection{Parallel Scalability}
	
	Finally, we investigate the parallel performance of both solvers.
	We consider again the shear test, but similar results can be expected for other scenarios as well.
	For this test, we take the mesh at $\ell = 9$, and consider  $p=1,2,3$ and $4$ leading to problem sizes 
	of {$7.9 \cdot 10^5$, $3.2 \cdot 10^6$, $7.1 \cdot 10^6$, and $1.3 \cdot 10^7$ dofs, respectively,}
	see \cref{tab:slit:dofs:p}.

	We notice that, for low degrees, the AMG approach is faster than our MF-GMG solver, 
	which we have already observed in \cref{fig:compare:time}.
	However, the solution time for the GMG solver grows 
	{slower than that for AMG solver with increasing $p$.}
	
	The parallel scaling behavior is similar for both preconditioners.
	For large enough problem sizes, strong scaling is close to perfect.
	As soon as the local size 
	{gets too small for each CPU, the}
	communication overhead starts to become noticeable.
	Both show good scaling down to approximately $25k$ dofs, although the AMG solver seems to perform slightly better.
	This is due to the bad scaling behavior on the coarse grids of the GMG hierarchy, which only contain very few cells.

	Both approaches do not show perfect weak scaling behavior.
	However, due to the complexity of the PDE considered here, these are more than satisfying results.
	
	We mentioned earlier, that memory consumption is a huge advantage of matrix-free methods.
	In fact, the matrix-based AMG simulation for $p=4$ exceeded our available RAM ($128$ GB) on a single core.
	
	\begin{figure}[H]
		\centering
		\includegraphics{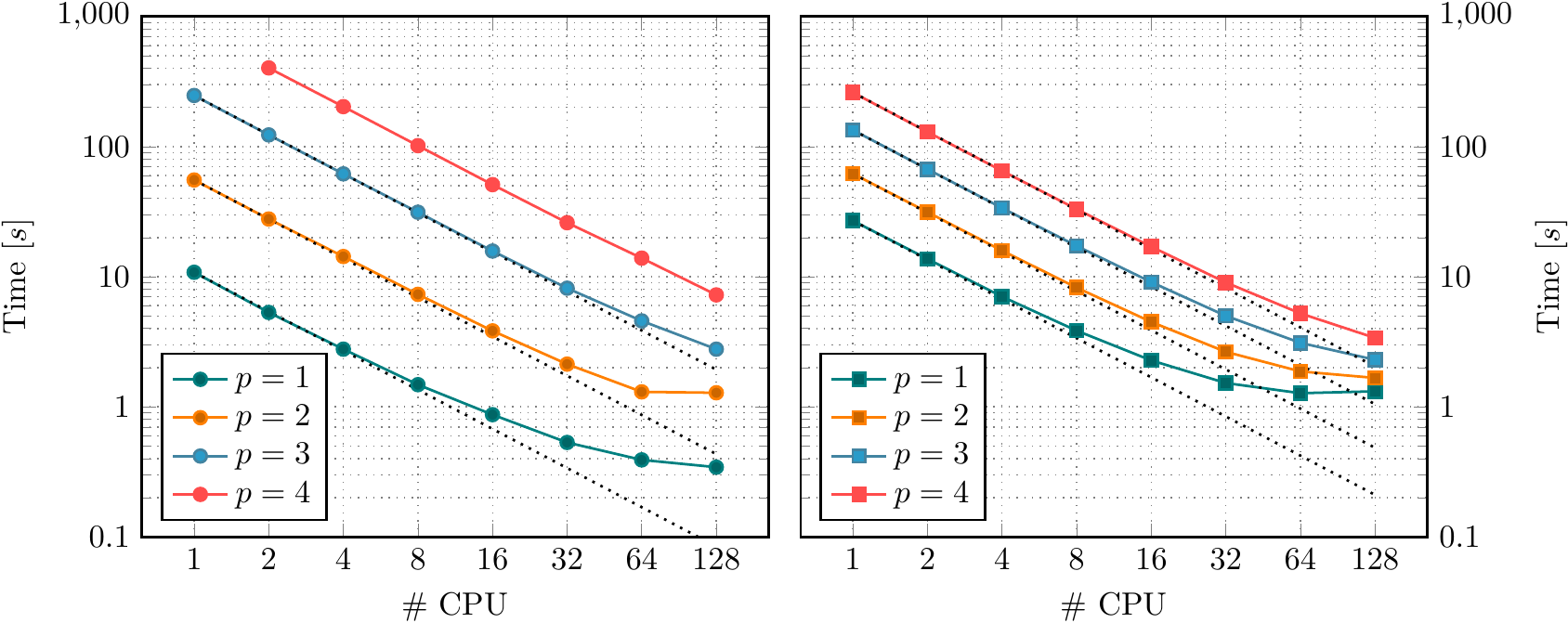}
		\caption{
			Parallel scalability of (\textbf{a}) AMG and (\textbf{b}) GMG for varying $p$.
			Average time spent for the solution of a linear system is plotted.
		}
		\label{fig:compare:scaling}
	\end{figure}

\section{Conclusions}

	We presented and compared two approaches to solve the linear systems arising in the PFF problem.
	First, we considered our matrix-free based geometric multigrid implementation from \cite{JoLaWi19}.
	Compared to our previous work, we extended our solver to handle higher polynomial degrees, which is a huge improvement for the matrix-free approach.
	
	We compared 
	{our}
	MF-GMG solver to the AMG-based approach 
	{from}
	\cite{HeWi18}.
	This method is comparably easy to implement
	{since}
	the AMG solvers 
        by MueLu (multigrid library in Trilinos) \cite{BeGlHu19}
	are almost black-box algorithms.
	These only require the sparse matrix, and compute the multigrid hierarchies by themselves.
	Implementing the MF-GMG approach is a lot more involved
	{since}
	we have to define the operators on each level.
	This can be particularly challenging in the presence of nonlinear terms and varying constraints (active-set), as it is the case in PFF.

	The matrix-free approach really excels 
	{at}
	high-order polynomial elements. 
	It requires less storage per dof as we increase the {polynomial} order $p$, whereas the sparse-matrix becomes more and more costly.
	This is also reflected in the performance of the linear solver, i.e., the time spent {for} solving the equations using AMG increases faster than using GMG during $p$-refinement.
	Nonetheless, the AMG approach is slightly faster for degrees less than $3$.
	Furthermore, the AMG solver requires a lot less implementational effort.
	The situation shifts more in favor of the matrix-free approach in $3$ dimensions, as storage costs tend to be even more limiting.

	So far, we only considered parallelization using SIMD instructions and distributed computing via MPI.
	Fine-grained parallelization within a node using threads (e.g., by means of OpenMP, TBB, std::thread) could improve {the} performance even more.
	An entirely different programming model is given by Graphic Processing Units (GPUs) using CUDA.
	GPUs are able to run several thousand threads in parallel, giving them a huge boost in performance for suitable application.
	Currently, an extension of the matrix-free framework in deal.II using GPUs is under development.

\section{Acknowledgments}
	
	This work has been supported by the Austrian Science Fund (FWF) grant P29181 `Goal-Oriented Error Control for Phase-Field Fracture Coupled to Multiphysics Problems'.

	\bibliographystyle{abbrv}
	\bibliography{literature}
	
\end{document}